\pdfoutput=1  

\documentclass{birkjour}

\usepackage{amssymb}
\usepackage{graphicx}
\usepackage[T1]{fontenc}
\usepackage[utf8]{inputenc}
\usepackage{url}
\usepackage{csquotes}
\usepackage{tikz}  

\makeatletter
\long\def\@firstoffive#1#2#3#4#5{#1}
\long\def\@secondoffive#1#2#3#4#5{#2}
\def\@newref#1{\expandafter\@setref\csname r@#1\endcsname\@firstoffive{#1}}
\def\@newpageref#1{\expandafter\@setref\csname r@#1\endcsname\@secondoffive{#1}}
\makeatother

\renewcommand{\subjclassname}{Mathematics Subject Classification (2020)}

\makeatletter
\renewcommand\paragraph{\@startsection{paragraph}{4}{\z@}%
  {-1.4ex \@plus -.4ex \@minus -.2ex}{-.5em}{\normalfont\bfseries}}
\makeatother

\newtheorem{proposition}{Proposition}

\newif\ifarxiv
\arxivtrue

\newcommand{\Pos}{\mathsf{P}}
\newcommand{\God}{\mathsf{G}}
\newcommand{\NE}{\mathsf{NE}}
\newcommand{\HF}{\mathit{HF}}
\newcommand{\wa}{\forall^{E}}
\newcommand{\we}{\exists^{E}}

\newcommand{\isabreakok}{\renewcommand{\_}{\textunderscore\allowbreak}}
\DeclareRobustCommand{\thy}[1]{{\isabreakok\textsf{#1}}}
\DeclareRobustCommand{\isaref}[2]{{\isabreakok\textsf{#1}\textup{.}\allowbreak\textsf{#2}}}

\begin{document}

\title[Demodalising G\"odel's God?]{Demodalising G\"odel's God?
       Ultrafilters, Rigidity and Modal Collapse Revisited}

\author[C.~Benzm\"uller]{Christoph Benzm\"uller}

\address{Otto-Friedrich-Universit\"at Bamberg, Bamberg, Germany, and
         Freie Universit\"at Berlin, Berlin, Germany}

\email{christoph.benzmueller@uni-bamberg.de}

\thanks{A comment on Odifreddi and Gomes, \emph{Logica Universalis}
        \textbf{20} (2026), 97--112, \texttt{doi:10.1007/s11787-026-00409-6}. Declaration:
        generative AI (Anthropic's Claude) supported the preparation of
        this paper, in particular, the machine-checked companion
        formalisation; the author has reviewed and verified all content,
        including every formal claim, and takes full responsibility for
        the publication.}

\subjclass{Primary 03A05; Secondary 03B45, 03B15, 68T15}

\keywords{Ontological Argument, Modal Ultrafilter, Modal Collapse,
          Computational Metaphysics, Higher-Order Modal Logic,
          Interactive and Automated Theorem Proving}

\begin{abstract}
G\"odel's ontological argument is \emph{modal} by design. Possibility and
necessity are among its central concerns, and connected interpretation
decisions --- extensions or intensions for the positive properties,
constant or varying quantifier domains, the choice of modal logic, etc.\ ---
subtly shape the argument and, in particular, whether the modal collapse is
implied or not. It is tempting to abstract the modality away and to read
the positive properties as merely forming an ultrafilter --- which,
however, may replace G\"odel's argument with a perhaps more manageable
alternative that may have little to do with G\"odel's original
intention.
In this spirit, Odifreddi
and Gomes argue that modal collapse is an intrinsic structural feature of
any theory that characterises positive properties as an ultrafilter and
identifies God as its principal generator. This comment shows that the
structural claim is false. There is an explicit $S5$ model, with constant
domains, in which positivity is at every world a modal ultrafilter
principally generated by a unique and necessarily existing Godlike
individual, and in which the modal collapse nevertheless fails --- already
for an atomic formula. What drives the collapse is not the ultrafilter
structure but the \emph{rigidity} of positivity --- itself a thoroughly
modal condition: given only that positivity is a modal ultrafilter, modal
collapse and rigidity are \emph{equivalent}. Rigidity is thus exactly what
the collapse adds to the ultrafilter structure, and the disagreement is
whether it is to be counted part of that structure at all. Two further claims of Odifreddi and Gomes are clarified: one on the
derivability of G\"odel's Theorem~IV, one on the extensionality of
positivity. All results are machine-checked in Isabelle/HOL, which
distinguishes kernel-checked proofs from model-finder output.
\end{abstract}

\maketitle

\section{Introduction}
\label{sec:intro}

The fact that G\"odel formulated his ontological argument in modal logic was
a deliberate choice, not a matter of notational convenience. The use of the
modalities possibility and necessity is one of his central concerns, and the
related interpretive decisions are complex, including, for example: whether
positive properties are extensions or intensions, whether quantifiers range
over constant or variable domains, whether full or only restricted
comprehension is available for properties, or which modal logic --- S5, S4,
KB, K --- is assumed (cf.\ \cite{Kovac2003,BenzmuellerScott2025,VestrucciBenzmueller2024}). These decisions subtly
influence the argument itself, its persuasiveness and its further
implications. It may therefore be tempting --- but not entirely harmless ---
to disregard the modalities and simply regard the positive properties as
forming an ultrafilter. This step risks replacing G\"odel's argument with a
more manageable substitute; since the disagreement discussed below stems
directly from this risk, the risk is examined first.

Two interesting aspects of G\"odel's argument have long been the subject of
debate. The first is the \emph{modal collapse}: from G\"odel's axioms one
can derive $\varphi \supset \Box\varphi$\footnote{Odifreddi and Gomes
\cite{OdifreddiGomes2026} state the collapse as the provability of
$\varphi \leftrightarrow \Box\varphi$ in their Section~5, and as
$\varphi \supset \Box\varphi$ in their Section~2.3. I use the
implication throughout:
the two coincide in any logic containing the reflexivity principle~T, and a
countermodel to the implication is a fortiori a countermodel to the
biconditional (\isaref{HOML}{MC\_iff\_T}, \isaref{HOML}{MC\_bicond\_imp}).} --- whatever is the case is necessarily the case --- as
noted already by Sobel \cite{Sobel1987}. The second is that G\"odel's
\emph{positive properties} --- his only primitive concept; colloquially, the
properties a perfect being should possess --- behave like an
\emph{ultrafilter}, which may loosely be described as a maximally coherent,
self-contained, rational structure (the formal definition follows in
Section~\ref{sec:ultrafilters}). The ultrafilter connection has been known
for quite some time \cite[pp.~364--365]{Hazen1998}, \cite{Odifreddi2000}; what is comparatively
recent, however, is that both aspects have been made precise and checked
using interactive and automated theorem proving technology
\cite{BenzmuellerWP2014,BenzmuellerFuenmayor2020,Benzmueller2020KR}.
Odifreddi and Gomes \cite{OdifreddiGomes2026} take the ultrafilter reading
as their starting point and put it to a new use: they connect it to
Friedman's \cite{Friedman2012} relative consistency results for set theory,
which rest on a \emph{non-principal} ultrafilter --- one not generated by
any single object. That part of their paper opens genuinely new ground and
is not at issue here.

At the close of their Section~2.3, Odifreddi and Gomes
\cite[\S2.3]{OdifreddiGomes2026} adopt the following simplification (the
reference numbering within the quotation has been adjusted to the present
bibliography):
\begin{quote}
As pointed out by several authors --- for a common reference, see Hazen \cite{Hazen1998} G\"odel's ontological argument forms an ultrafilter
over the set of positive properties. Consequently, the modal vocabulary can
be dispensed with in the analysis of its deep structure.
\end{quote}
As a practical approach to the consistency-theoretic goals of their
Section~4, this may make sense, and nothing here calls it into
question. As a statement about \emph{G\"odel's} argument, however, the
simplification has more far-reaching consequences than it may seem at first
glance. Whether the (modal) ultrafilter of positive properties is defined
via property \emph{extensions} or via property
\emph{intensions} is not a mere side issue: it determines the type
of the modal ultrafilter, and it bears directly on the modal collapse. Restricting the
ultrafilter condition to extensions, as in Fitting's variant
\cite{Fitting2002} --- a reading Fitting introduced precisely to
\emph{avoid} the collapse --- does avoid it (machine-checked in
\cite{BenzmuellerFuenmayor2020}), while the intensional reading does not by
itself force it (Section~\ref{sec:counterexamples}).\footnote{The quoted
passage appeals to Hazen \cite{Hazen1998}; the extensional reading Odifreddi
and Gomes themselves work with is Friedman's \cite{Friedman2012}, adopted in
their Section~4.} Setting the modality aside is therefore a change of system
--- and one would like to be sure that the theory shown to collapse is still
G\"odel's.

In their Section~5.1, Odifreddi and Gomes in contrast hold that modal
collapse $\varphi \supset \Box\varphi$ is
\begin{quote}
not an accidental by-product of an overly strong modal logic; it is an
intrinsic feature of any argument that characterises positive properties as
an ultrafilter and then identifies God as the (principal) generator of that
ultrafilter
\end{quote}
--- and, correlatively, that avoiding the collapse by weakening
the axioms costs the principality of the ultrafilter and hence the absolute
maximality of the divine being (their Sections~5.1 and~5.2). To be
fair, their surrounding prose is partly more guarded --- the ultrafilter
reading \enquote{suggests} that the collapse \enquote{may instead reflect}
a deeper structural feature of the positivity conditions --- but the
passage quoted above is their most explicit formulation, it is the one
their Sections~5.2 and~6 build on, and it is therefore the claim examined
here.

This claim faces a dilemma. Read modally --- as a statement about theories
in higher-order modal logic in which positivity forms a (modal) ultrafilter
whose principal generator is Godlikeness --- the claim is false:
Section~\ref{sec:counterexamples} provides verified counterexamples. Read
extensionally --- as a claim about the demodalised, set-theoretic surrogate
--- there is no modal collapse left to assert; what collapses there is the
demodalised surrogate rather than G\"odel's argument, and
Section~\ref{sec:ext} specifies where the transition takes place. (If \emph{intrinsic} is instead
meant explanatorily, it escapes formal refutation --- but then also loses
the argumentative weight of their Sections~5.2 and~6.)

\paragraph{Companion formalisation.}
All claims below are supported by Isabelle/HOL \cite{Isabelle} theory
files, based on the embedding of \cite{BenzmuellerPaulson2013} and
available as ancillary files with the arXiv version of this
article (\texttt{arXiv:2608.07578}).\footnote{The theory files provided are \thy{HOML}, \thy{MFilter},
\thy{BaseDefs}, \thy{ScottVariant}, \thy{AndersonVariant},
\thy{FittingVariant}, \thy{UFilterVariant}, \thy{UFilterS5},
\thy{UFilterOnly}, \thy{SimpleVariantHF},
\thy{GoedelAx1Gen}, \thy{GoedelInconsistency}, \thy{ExtensionalityTests},
\thy{Th4Underivability} and \thy{AltEntailmentVariant}.} To keep the main
text self-contained, the verification pointers \isaref{Theory}{lemma} that
name the checking statement are relegated to footnotes. Countermodels are
found automatically, by the model finder Nitpick \cite{Nitpick}, which is
much of what makes this setting productive to work in. A model so found can
then be written out explicitly, as an interpretation of the types and
constants of the embedding, and the countermodel claim \emph{proved} about
it, after which the result no longer depends on the finder. Both steps take
place within Isabelle/HOL, and both are used below:
Section~\ref{sec:ufvariant} displays a model that began as a Nitpick
countermodel and is now a theorem. Section~\ref{sec:repro} records which claim
rests on which.

\section{Two Notions of Modal Ultrafilter}
\label{sec:ultrafilters}

Everything below rests on one distinction. The property \emph{being the
tallest person in the room} may, as a concept, pick out different
individuals in different circumstances: the concept is the property's
\emph{intension}, the set of individuals it picks out in one fixed
circumstance its \emph{extension}. A property is \emph{rigid} if its
intension determines the same extension in every circumstance. The
following types, as used in
\cite{BenzmuellerFuenmayor2020,BenzmuellerPaulson2013}, make this precise.
Let $i$ be the type of possible worlds, and $e$ the type of individuals.
A \emph{proposition}, of type $\sigma := i \Rightarrow \mathit{bool}$, may be true at some
worlds and false at others; a property \emph{extension}, of type
$\delta := e \Rightarrow \mathit{bool}$,
is a set of individuals; a property
\emph{intension}, of type $\gamma := e \Rightarrow \sigma$, is a world-dependent
concept, determining an extension at each world. A proposition $\psi$ of
type $\sigma$ is
\emph{valid}, written $\lfloor\psi\rfloor$, if it
holds at every world; these definitions are presented in the Isabelle file
\thy{HOML}. Sets are represented throughout by
predicates (characteristic functions), a deliberate difference from the
set-theoretic ultrafilter of \cite{OdifreddiGomes2026} that allows the
distinctions below to be expressed at all.
\emph{Possibilist} quantifiers range over all individuals; \emph{actualist}
quantifiers, written $\wa$ and $\we$, only over the individuals
\emph{existing} at the world of evaluation.

\paragraph{The classical notion and its modal counterpart.}
The filter and ultrafilter of \cite[Def.~3.1, conditions (4)--(5)]{OdifreddiGomes2026}
are the classical global notions on a fixed index set~$I$: a family of
subsets of $I$ that contains $I$, omits $\emptyset$, and is closed under
supersets and finite intersections; an ultrafilter contains in addition
exactly one of each set and its complement. Taking $I$ to be the individuals
and reading positivity extensionally yields the extensional notion, written
$\delta$ below. A modal filter is that same notion with a single change: the conditions are imposed
\emph{afresh at each world}, inclusion being evaluated at the world of
evaluation over the individuals existing there (cf.\
\cite[Def.~1.5]{Kovac2003}). Written out, a world-indexed family $\phi$ of
properties is a \emph{modal filter} if at every world
\[
\begin{aligned}
  &\lambda x.\top \in \phi, \qquad \lambda x.\bot \notin \phi, \\[2pt]
  &(X \in \phi \wedge X \sqsubseteq Y) \supset Y \in \phi, \qquad
   (X \in \phi \wedge Y \in \phi) \supset X \sqcap Y \in \phi,
\end{aligned}
\]
where $X \sqsubseteq Y := \wa z.(X z \supset Y z)$ and
$(X \sqcap Y)\,z := X z \wedge Y z$, and a \emph{modal ultrafilter} if in
addition $X \in \phi \vee {\sim}X \in \phi$ at every world, for every
$X$.\footnote{\thy{MFilter}.}

That this is the natural modal counterpart of the classical notion can be
put precisely: a world-indexed family is a modal (ultra)filter exactly when, at each world, membership depends only on the local extension and the local extensions of its members form a classical (ultra)filter over the individuals existing there. The classical notion is
thus recovered at every single world, while the family as a whole is left
free to vary from world to world. Nothing in the classical definition settles
that variation either way, because over a fixed index set there is nothing to
vary; insisting that it not vary along accessibility is therefore a further
and substantive demand --- G\"odel's rigidity axiom \textbf{A4}
(Section~\ref{sec:a4}).

Two kinds of modal (ultra)filter thus arise \cite{BenzmuellerFuenmayor2020},
according to the type over which positivity is taken: a
\emph{$\gamma$-ultrafilter} weighs intensions --- world-dependent concepts
--- and is classical only through the local extensions of its members, as
just described, while a \emph{$\delta$-ultrafilter} weighs extensions ---
plain sets of individuals --- and is at each world a classical ultrafilter as
it stands. The former is the type G\"odel's own $\Pos$ has, and the
distinction is one that a single set-theoretic notion cannot express.

\paragraph{Principality.} The \emph{principal filter} generated by a
property $X$, written $\HF_X$ (for \emph{Hauptfilter}), is the family of all
supersets of $X$, again read afresh at each world. Classical usage is worth
keeping in view: the family of all supersets of a nonempty set is a filter,
principal by definition, whatever the generator, and it is an
\emph{ultra}filter precisely when the generator is a singleton. Over a fixed domain, then,
\enquote{principal ultrafilter generated by God} does carry uniqueness, and
hence monotheism, whereas \enquote{principal filter generated by God} by
itself does not; Section~\ref{sec:hf} keeps the two apart. A further
distinction has no classical analogue and matters here: positivity may be
generated at each world by some Godlike individual \emph{of that world}, or
by one and the same individual at \emph{every} world. The model of
Section~\ref{sec:ufvariant} realises the stronger, world-invariant form.

Finally, write $\Pos'$ for the derived test on intensions defined as
follows: at world $w$, $\Pos' X$ holds iff the \emph{rigidification} of
$X$'s extension at $w$ --- the constant concept that picks out, at every
world, exactly the individuals falling under $X$ at $w$ --- is positive at
$w$. Table~\ref{tab:variants} summarises the results.

\begin{table}[t]
\caption{Ultrafilter status of positivity across the classical variants,
after \cite{BenzmuellerFuenmayor2020}. The type column resolves what
\enquote{the ultrafilter} refers to in each row: in Fitting's variant $\Pos$
itself ranges over extensions, so the comparison $\Pos = \Pos'$ does not
apply and the $\delta$-notion takes the place of the $\gamma$-notion. The
principality column concerns the ultrafilter that the row supplies --- $\Pos$
for Scott, the derived $\Pos'$ for Anderson and Fitting --- and records
whether it is, at each world, the principal filter generated by $\God$.
\enquote{Collapse} records whether $\varphi \supset \Box\varphi$ is
derivable. A \enquote{yes} is a derivability over the embedding's actualist quantifiers --- varying domains, hence a fortiori constant ones --- in the weakest base in which it was checked: symmetry ($KB$) for the Scott row, $S5$ as axiomatised in the ancillary files for the Anderson and Fitting rows; the countermodels behind each \enquote{no} satisfy the $S5$ frame conditions and have constant domains.}
\label{tab:variants}
\centering
\footnotesize
\setlength{\tabcolsep}{4pt}
\begin{tabular}{llccccc}
\hline\noalign{\smallskip}
Variant & type of $\Pos$ & $\Pos$ an & $\Pos'$ a & $\Pos = \Pos'$ &
principal, & Collapse \\
 & & ultrafilter & $\gamma$-ultrafilter & & gen.\ by $\God$ & \\
\noalign{\smallskip}\hline\noalign{\smallskip}
Scott \cite{Scott1972}       & $\gamma$ & yes & yes & yes & yes & \textbf{yes} \\
Anderson \cite{Anderson1990} & $\gamma$ & \textbf{no} & yes & no & yes & \textbf{no} \\
Fitting \cite{Fitting2002}   & $\delta$ & yes ($\delta$) & yes & --- & yes & \textbf{no} \\
\noalign{\smallskip}\hline
\end{tabular}
\end{table}

The results of the table are due to
\cite{BenzmuellerFuenmayor2020,FuenmayorBenzmueller2017,FuenmayorBenzmuellerAFP2017};
all entries are replayed in the ancillary files.\footnote{Replayed
in \thy{ScottVariant} (\isaref{ScottVariant}{U1}, \isaref{ScottVariant}{U2},
\isaref{ScottVariant}{U3}, \isaref{ScottVariant}{MC}),
\thy{AndersonVariant} (\isaref{AndersonVariant}{U1\_fails},
\isaref{AndersonVariant}{U2}, \isaref{AndersonVariant}{U3\_fails},
\isaref{AndersonVariant}{MC}) and \thy{FittingVariant}
(\isaref{FittingVariant}{U1}, \isaref{FittingVariant}{MC}).} In
Fitting's variant, where positivity itself ranges over extensions, $\Pos'$
is likewise definable --- $\Pos' X$ holds at $w$ iff the extension of $X$
at $w$ is positive there --- and it, too, forms a
$\gamma$-ultrafilter:\footnote{\isaref{FittingVariant}{U2}.} the ultrafilter structure of the
$\delta$-row is thus visible at type $\gamma$ as well. An
ultrafilter structure hence survives in the collapse-free variants,
migrating from the intensions to the (rigidified) extensions.

The principality column answers, in passing, a claim that
\cite[\S\S5.1, 5.2]{OdifreddiGomes2026} makes about the collapse-free
emendations: that they must render the ultrafilter non-principal and thereby
forfeit the absolute maximality of God. They need not. In Anderson's and in Fitting's variant the relevant ultrafilter, $\Pos'$, is at every world principally generated by Godlikeness --- a property belongs to it exactly when every Godlike individual has it --- just as $\Pos$ is in Scott's, and the collapse nevertheless fails.\footnote{\isaref{AndersonVariant}{PrincipalA},
\isaref{FittingVariant}{PrincipalF}; for comparison
\isaref{ScottVariant}{Principal} and \isaref{ScottVariant}{Principal\_HF}.}
Principality is therefore not the price of avoiding the collapse.

And when
\cite[\S3.1]{OdifreddiGomes2026} speaks of \enquote{the ultrafilter of all
positive properties}, the referent shifts between a $\gamma$-ultrafilter, a
$\delta$-ultrafilter, and an ordinary set-theoretic ultrafilter over the
individuals --- the difference being exactly what the disputed claim turns
on.

\section{Ultrafilters Without Collapse; Principal Filters Without Collapse}
\label{sec:counterexamples}

\subsection{Taking the ultrafilter property as an axiom}
\label{sec:ufvariant}

Since G\"odel's premises \emph{entail} that $\Pos$ is a modal
ultrafilter,\footnote{\isaref{ScottVariant}{U1}.} one may reverse the order of business and
postulate outright, as axiom \textbf{U1}, that $\Pos$ is a modal
ultrafilter, obtaining one of the deliberately simplified theories of
\cite{Benzmueller2020KR}. Retaining only \textbf{A2} (a property necessarily
entailed by a positive property is positive), \textbf{A3} (any conjunction
of positive properties is positive) and the definition
$\God x \equiv \forall Y.(\Pos Y \supset Y x)$, one
obtains:\footnote{Three labelling systems are in play --- Scott's, G\"odel's
and that of \cite{OdifreddiGomes2026} --- and they do not agree.
Table~\ref{tab:concordance} is a concordance. On the textual history of
Scott's replacement of the conjunction axiom by $\Pos(\God)$ see
\cite{BenzmuellerScott2025}.}

\begin{table}[t]
\caption{Concordance of the three labelling systems. The column
\enquote{Here} follows Scott's numbering with one deliberate difference, marked
$(\dagger)$: where Scott postulates $\Pos(\God)$ outright as his
\textbf{A3}, the \textbf{A3} used in
Section~\ref{sec:counterexamples} is the stronger conjunction principle of
\cite{Benzmueller2020KR}, from which $\Pos(\God)$ is derived. G\"odel's
\textbf{Ax1} is binary in the scriptum, and extended in his own marginal
footnote to any number of conjuncts; \cite{OdifreddiGomes2026} states
Axiom~I with that footnote and discusses the distinction. Dashes mark principles for which the system in question has no label; in the column \enquote{Here}, \textbf{A1}, \textbf{A4} and \textbf{A5} are the labels under which Section~\ref{sec:counterexamples} discusses the principles it drops.}
\label{tab:concordance}
\centering
\footnotesize
\setlength{\tabcolsep}{4.5pt}
\begin{tabular}{lcccc}
\hline\noalign{\smallskip}
Principle & Here & Scott & G\"odel & O\&G \\
\noalign{\smallskip}\hline\noalign{\smallskip}
binary conjunction closure
  & --- & --- & \textbf{Ax1} & I \\
arbitrary conjunction closure
  & \textbf{A3}$^{\dagger}$ & --- & \textbf{Ax1Gen} & I (fn.) \\
$\Pos(\God)$
  & --- & \textbf{A3} & derived & Ax.\ 3.6 \\
exactly one of $X$, ${\sim}X$ positive
  & \textbf{A1} & \textbf{A1} & \textbf{Ax2a} & IIA \\
rigidity: $\Pos X \supset \Box(\Pos X)$
  & \textbf{A4} & \textbf{A4} & \textbf{Ax2b} & IIB \\
closure under necessary entailment
  & \textbf{A2} & \textbf{A2} & \textbf{Ax4} & IV \\
$\Pos(\NE)$
  & \textbf{A5} & \textbf{A5} & \textbf{Ax3} & III \\
\noalign{\smallskip}\hline
\end{tabular}
\end{table}

\begin{proposition}[after \cite{Benzmueller2020KR}, \S5; \thy{UFilterVariant},
\thy{UFilterS5}]
\label{prop:ufvariant}
From \textbf{U1}, \textbf{A2}, \textbf{A3} and $\mathbf{df.}\God$ the
theorem $\lfloor \Box \we x.\,\God x \rfloor$ --- at every world, some
individual existing there is Godlike --- is derivable already in modal
logic~$K$, the weakest normal modal logic. G\"odel's remaining axioms may be
dropped: \textbf{A1} --- that of $X$ and ${\sim}X$ exactly one is positive,
G\"odel's \textbf{Ax2a} --- stays derivable from \textbf{U1} alone, while
\textbf{A4} and \textbf{A5} do not. The theory is consistent, and modal collapse fails in an
explicit $S5$ model with constant domains.\footnote{In order of appearance:
\isaref{UFilterVariant}{T6}; \isaref{UFilterVariant}{A1\_implied};
\isaref{UFilterVariant}{A4\_not\_implied}, \isaref{UFilterVariant}{A5\_not\_implied};
\isaref{UFilterS5}{MC\_atom\_fails}, \isaref{UFilterS5}{MC\_fails}.}
\end{proposition}

It is worth being explicit about what is retained, what is dropped, and what
is thereby refuted. Retained are \textbf{U1} --- positivity is a modal
ultrafilter, the conditions imposed afresh at each world --- together with
\textbf{A2}, \textbf{A3} and the definition
$\God x \equiv \forall Y.(\Pos Y \supset Y x)$, so that Godlikeness remains
\emph{maximal}: the having of \emph{all} positive properties. Dropped are
\textbf{A1} (which returns as a theorem), \textbf{A4}, \textbf{A5} and the
definitions of essence and of necessary existence. The claim refuted is
accordingly the following: that a theory in which positivity is at each world
a modal ultrafilter principally generated by a Godlike individual must
validate $\varphi \supset \Box\varphi$. It is \emph{not} claimed that the
collapse fails in the full G\"odel--Scott system, where it is
derivable,\footnote{\isaref{ScottVariant}{MC}.} nor that every G\"odelian
premise survives here.

Two facts about the hypotheses themselves should be on record. First,
\textbf{A2} does no work in the derivation of
$\lfloor \Box \we x.\,\God x \rfloor$: $\Pos(\God)$ follows from
\textbf{A3}; a Godlike individual then exists at every world by \textbf{U1}
together with $\Pos(\God)$, since were the extension of $\God$ empty among the
existents the empty property would be in the filter by superset closure from
$\Pos(\God)$; and
necessary existence follows by necessitation. \textbf{A2} is kept because the
theory of \cite{Benzmueller2020KR} has it, and because the countermodel is
thereby shown to satisfy it too. Second, \textbf{A3} cannot be dropped. The
infinite model of Section~\ref{sec:thmIV} satisfies \textbf{U1} and
\textbf{A2} --- indeed also rigidity --- and has no Godlike individual at
all; what it lacks is precisely closure under \emph{arbitrary}
conjunctions.\footnote{\isaref{UFilterVariant}{T2},
\isaref{UFilterVariant}{ExistsG}, \isaref{UFilterVariant}{T6} for the
derivation; \isaref{Th4Underivability}{U1}, \isaref{Th4Underivability}{Ax4}
and \isaref{Th4Underivability}{Th4\_fails} for the model that separates
\textbf{A3} from the rest.}

The countermodel is small enough to state in full. Nitpick finds one
automatically, and that is how it was first obtained; what follows is that
model written out, with its accessibility relation strengthened to the
\emph{universal} frame, since G\"odel, Scott and \cite{OdifreddiGomes2026} all
work in $S5$ and, for a negative claim, $K$ is the weakest setting rather than
the strongest. Written out, it is no longer the finder's output but an
interpretation about which everything asserted below is proved. Take two worlds $i_1$, $i_2$ and a single
individual $e_1$; let every world access every world, so that accessibility
is an equivalence relation and the model is an $S5$ model, a fortiori one of
$K$, $KB$ and $S4$; let $e_1$ exist at both worlds, so that the domains are
constant; let $e_1$ be Godlike at both worlds; and let the positive
properties at a world be exactly the properties $e_1$ has there,
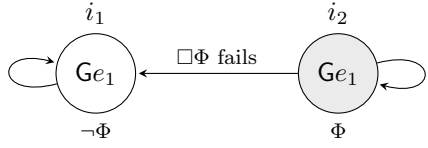
\[
  \Pos(X,w) \quad\text{iff}\quad X(e_1,w).
\]
\begin{center}
\begin{tikzpicture}[->,>=stealth,shorten >=1pt,
  world/.style={draw,circle,minimum size=1.05cm,inner sep=1pt},
  phi/.style={font=\footnotesize}]
  \node[world,label=above:{$i_1$},
        label={[phi]below:{$\lnot q$}}] (i1) at (0,0) {$\God e_1$};
  \node[world,fill=black!8,label=above:{$i_2$},
        label={[phi]below:{$q$}}] (i2) at (3.4,0) {$\God e_1$};
  \path (i1) edge [loop left]  (i1)
        (i2) edge [loop right] (i2)
        (i1) edge [bend left=20]  (i2)
        (i2) edge [bend left=20]  (i1);
\end{tikzpicture}
\end{center}
(Accessibility is universal, hence reflexive, symmetric and transitive. The
shaded world marks the extension of the atom $q := \{i_2\}$ discussed
below.)

Positivity so construed satisfies the ultrafilter conditions at both worlds:
the universal property is positive and the empty one is not, since $e_1$
falls under the former and not under the latter; closure under supersets and
under finite intersections is immediate; and maximality holds because of each
property and its complement $e_1$ has exactly one. So \textbf{U1} holds.
\textbf{A2} holds because each world accesses itself: if $X$ is positive at
$w$ and $X$ necessarily entails $Y$, then in particular $X \sqsubseteq Y$ at
$w$, so $Y$ holds of $e_1$ there. \textbf{A3} --- every, possibly infinite,
conjunction of positive properties is positive --- holds for the same reason:
if $X$ is the conjunction of a family all of whose members are positive, then
every member holds of $e_1$ at $w$, hence so does $X$. Finally $\God e_1$
holds at both worlds, since a property positive at $w$ holds of $e_1$ at $w$
by the very definition of $\Pos$.

Yet contingency survives. Extend the language by a propositional constant
$q$, interpreted as the singleton $\{i_2\}$: true at $i_2$, false at $i_1$.
($q$ is an ordinary atom, interpreted by a set of worlds as any atom is; that
the set is a singleton is a feature of this two-world model, not a nominal of
the language.)
At $i_2$ the atom $q$ holds while $\Box q$ fails, $q$ being false at the
accessible world $i_1$; so the collapse instance $q \supset \Box q$ is
invalidated, and the collapse is refuted already in its schematic reading,
for an \emph{atomic} formula. Read instead as the quantified statement
$\lfloor \forall \Phi.\,\Phi \supset \Box\Phi \rfloor$ --- the reading on
which $\Phi$ ranges over all propositions of the model --- the collapse fails
at \emph{every} world.\footnote{\isaref{UFilterS5}{MC\_atom\_fails},
\isaref{UFilterS5}{MC\_fails}; the frame conditions are
\isaref{UFilterS5}{R\_refl}, \isaref{UFilterS5}{R\_symm},
\isaref{UFilterS5}{R\_trans} and constancy of the domains
\isaref{UFilterS5}{constant\_domain}.} The ultrafilter property of positivity
is compatible with contingency.

Note, moreover, that the \emph{full} hypothesis of the disputed claim is
realised in this model. Positivity is at each world principally generated by
Godlikeness --- a property is positive iff every Godlike individual has it,
$\lfloor \forall X.\, \Pos X \leftrightarrow \forall x.(\God x \supset X x)
\rfloor$, equivalently $\Pos$ coincides with the principal filter
$\HF_\God$ --- the generator is the \emph{same} individual at every world,
not merely one per world; that individual is unique and exists necessarily;
and \textbf{U1} forces monotheism in general.\footnote{\isaref{UFilterS5}{Principal},
\isaref{UFilterS5}{Principal\_HF}, \isaref{UFilterS5}{MT},
\isaref{UFilterS5}{T6}; \isaref{UFilterVariant}{Principal},
\isaref{UFilterVariant}{MT\_U1}.} An ultrafilter of positive properties with
God as its principal generator, in $S5$, with constant domains --- and still
no collapse. Every assertion about the displayed model is checked by the
Isabelle kernel; no model finder is involved.

\subsection{What actually drives the collapse}
\label{sec:a4}

The same model yields a sharper diagnosis. G\"odel's rigidity axiom (the
first conjunct of Axiom~IIB of \cite{OdifreddiGomes2026}),
\textbf{A4}: $\forall X.\,(\Pos X \supset \Box(\Pos X))$, fails in it --- and
it fails at a witness of the object language. Take the constant property
$\lambda x.\,q$, whose extension is $\{e_1\}$ at $i_2$ and empty at $i_1$: it
is positive at $i_2$, while $\Box(\Pos(\lambda x.\,q))$ fails there, since
$\Pos(\lambda x.\,q)$ is false at the accessible world $i_1$. Like the
collapse, \textbf{A4} fails at every world of the
model.\footnote{\isaref{UFilterS5}{A4\_atom\_fails},
\isaref{UFilterS5}{A4\_fails}; likewise \textbf{A5},
\isaref{UFilterS5}{A5\_fails}.}

Adding \textbf{A4} back restores the collapse. One route runs through
G\"odel's own apparatus: with \textbf{A4}, Godlikeness is an essence of any
Godlike individual, so every property such an individual has --- in
particular, for a true proposition $Q$, the constant property
$\lambda y.\,Q$ --- is necessarily entailed by $\God$, and the necessary
existence of a Godlike individual then yields
$\Box Q$.\footnote{\isaref{UFilterVariant}{MC\_with\_A4}.} None of that
apparatus is needed, however: the modal ultrafilter conditions alone tie the
positivity of a constant property to the truth of its proposition, and
dispensing with the apparatus makes the diagnosis exact.

\begin{proposition}[\thy{UFilterOnly}]
\label{prop:equivalence}
Assume \textbf{U1} alone --- positivity is a modal ultrafilter --- in base
logic $K$: no further axiom, no frame condition, and neither essence nor
necessary existence defined. Then modal collapse and the rigidity of
positivity are equivalent,
\[
  \lfloor \forall \Phi.\,\Phi \supset \Box\Phi \rfloor
  \qquad\text{iff}\qquad
  \lfloor \forall X.\,\Pos X \supset \Box(\Pos X) \rfloor ,
\]
and \textbf{U1} implies neither.\footnote{\isaref{UFilterOnly}{MC\_iff\_A4},
from \isaref{UFilterOnly}{MC\_from\_A4} and
\isaref{BaseDefs}{A4\_from\_MC} (the latter postulates nothing, and is
placed accordingly); the two non-implications are witnessed by the explicit
model of Section~\ref{sec:ufvariant}, which satisfies \textbf{U1}
(\isaref{UFilterS5}{MC\_fails}, \isaref{UFilterS5}{A4\_fails}), and recorded
inside the theory by Nitpick (\isaref{UFilterOnly}{MC},
\isaref{UFilterOnly}{A4\_not\_implied}); \isaref{UFilterOnly}{U1\_A4\_consistent}
records that the hypothesis is satisfiable, so the equivalence is not
vacuous.}
\end{proposition}

Both directions are short. From left to right, instantiate the collapse at
$\Phi := \Pos X$; this uses no axiom at all. From right to left, note first
that a constant property $\lambda x.\,\varphi$ is positive at exactly those
worlds at which $\varphi$ is true. If $\lambda x.\,\varphi$ were positive at a world where $\varphi$
fails, it would have the empty extension there, and closure under supersets
would place the empty property in the filter, which \textbf{U1} forbids. If
$\varphi$ held at a world where $\lambda x.\,\varphi$ is not positive,
maximality would put the complement $\lambda x.\,\neg\varphi$ into the
filter, and by what was just shown $\neg\varphi$ would hold there --- a
contradiction. Now let $\varphi$ hold at $w$. Then $\lambda x.\,\varphi$ is
positive at $w$; by \textbf{A4} it is positive at every world accessible from
$w$; and there the first observation returns $\varphi$. Hence $\Box\varphi$
at $w$.

The right-to-left step reads inclusion locally while excluding the globally
empty property, and one might suspect that the equivalence trades on this
pairing. It does not: as Section~\ref{sec:ultrafilters} put it, a modal (ultra)filter is at each world a classical (ultra)filter of local extensions, and both features are what that requires. The equivalence therefore bears on the
ultrafilter reading itself, not on one formalisation of it.

What the equivalence does \emph{not} establish should be said as plainly. It
is a statement about theories containing \textbf{U1}: it does not say that
rigidity entails the collapse in every theory, nor that rigidity is necessary
for the collapse across theories. Kova\v{c} \cite{Kovac2003}, for instance,
blocks the general collapse without touching G\"odel's axioms, by
restricting comprehension: his $\lambda$-abstracts may contain no closed
subformula \cite[Def.~1.3, p.~569]{Kovac2003}, which removes the constant
property $\lambda y.\,\varphi$ on which Sobel's derivation --- and the proof
above --- turns; his $KB$ system keeps the axioms, with Scott's
$\Pos(\God)$, and still has the rigidity of positivity,
$\forall X.(\Pos X \leftrightarrow \Box\Pos X)$, as a theorem
\cite[Prop.~4.4, p.~575]{Kovac2003}. What the equivalence does establish is
that at the level of generality at which \cite{OdifreddiGomes2026} state the
structural thesis --- positivity an ultrafilter, God its principal generator
--- the dividing line is rigidity and nothing else. This fits
Table~\ref{tab:variants}: Fitting's emendation changes the type of $\Pos$,
and in Anderson's the $\gamma$-ultrafilter that survives, $\Pos'$, is not
rigid --- as, by the argument above, it cannot be, the collapse failing
there.\footnote{\isaref{AndersonVariant}{Prigid\_not\_rigid}; Anderson's
rigidity axiom concerns $\Pos$, which is not the ultrafilter
(\isaref{AndersonVariant}{U1\_fails}).}

It also locates the disagreement. On a natural reading of
\cite{OdifreddiGomes2026}, \enquote{the ultrafilter} is a single, fixed
ultrafilter, the same at every world, so that world-invariance is built into
the phrase. For G\"odel's $\Pos$, which ranges over intensions, that
invariance is rigidity, and on that reading their thesis is correct by
Proposition~\ref{prop:equivalence}; what divides us is then not a theorem
but whether rigidity is to be counted part of the \enquote{ultrafilter
structure}. (For extensions it is not correct even so: Fitting's $\Pos$ is
rigid and a $\delta$-ultrafilter, and the collapse fails,
Table~\ref{tab:variants}.) The present point is that rigidity does not come
for free: a notion of ultrafilter imposed afresh at each world --- the notion
one needs if the modal vocabulary is to be retained at all --- is satisfiable
without it, together with principal generation by a unique and necessarily
existing God (Section~\ref{sec:ufvariant}), and rigidity is exactly what the
collapse adds.

\subsection{Principality without collapse}
\label{sec:hf}

The second half of the disputed claim ties principality --- being generated
by a single object, here the Godlike being --- to the collapse. The modal
counterpart is the principal filter
$\HF_\God$ of all supersets of $\God$; G\"odel's theory, in
Scott's consistent form, implies that $\HF_\God$ is a modal
(ultra)filter.\footnote{\thy{MFilter}; \isaref{ScottVariant}{F\_HF}, \isaref{ScottVariant}{UF\_HF}.} Now assume
only axiom \textbf{F1}: $\HF_\God$ is a modal filter.

\begin{proposition}[\cite{Benzmueller2020KR}, \S7; \thy{SimpleVariantHF}]
\label{prop:hauptfilter}
From \textbf{F1} alone, the existence and the necessary existence of a
Godlike entity follow; the theory is consistent, and neither modal
collapse nor monotheism is implied.\footnote{\isaref{SimpleVariantHF}{T3'}, \isaref{SimpleVariantHF}{T6};
\isaref{SimpleVariantHF}{MC}, \isaref{SimpleVariantHF}{MT}; the countermodels
also survive the $S5$ frame conditions, \isaref{SimpleVariantHF}{MC\_S5},
\isaref{SimpleVariantHF}{MT\_S5}.}
\end{proposition}

What this proposition does and does not do should be stated plainly. In the
\textbf{F1} and \textbf{UF1} theories only $\HF_\God$ is constrained;
$\Pos$ itself is \emph{not} required to be an ultrafilter. These theories
therefore do not instantiate the hypothesis of
\cite{OdifreddiGomes2026} that the positive properties form an ultrafilter,
and they are not what refutes the structural thesis --- that is carried by
Proposition~\ref{prop:ufvariant} and the model of
Section~\ref{sec:ufvariant} alone. Their business is the second half of the
disputed claim, the one about principality and monotheism, and it is in that
role that they are reported here.

To be transparent: \textbf{F1}'s only substantive demand --- read afresh at each
world --- is that at every world some existing individual is Godlike;
existence is thus in effect postulated. This variant is accordingly the
most drastic of the simplifications considered here, and the absence of a
collapse in it is correspondingly unsurprising. Its interest lies in what does \emph{not}
follow. Strengthening \textbf{F1} to
\enquote{$\HF_\God$ is a modal \emph{ultrafilter}} does not change the
picture: the collapse still has a countermodel, while uniqueness among
existing Godlike entities is now forced --- as it must be, since a principal
\emph{ultra}filter is one whose generator is a singleton, which is why
\textbf{F1}, postulating a principal \emph{filter} only, leaves monotheism
open while \textbf{UF1} does not; and the ultrafilter condition of \textbf{U1} --- of $X$ and ${\sim}X$ one is positive --- likewise forces monotheism in the ultrafilter variant.\footnote{\isaref{SimpleVariantHF}{MC\_UF}, and under the $S5$ frame
conditions \isaref{SimpleVariantHF}{MC\_UF\_S5}; \isaref{SimpleVariantHF}{MT\_UF};
\isaref{UFilterVariant}{MT\_U1}.} Table~\ref{tab:simplified}
collects the verified facts.

\begin{table}[t]
\caption{The deliberately simplified theories at a glance. A \enquote{yes} is a derivability in base logic~$K$ over varying domains, a \enquote{no} a countermodel satisfying the $S5$ frame conditions with constant domains --- in each case the strongest claim of its kind, so that no entry depends on the choice of a weak modal base or of the domain semantics. In
the \textbf{U1} row the countermodel is the explicit $S5$ model of
Section~\ref{sec:ufvariant}, with constant domains and checked by the kernel; the other countermodels are found by Nitpick, with the $S5$ frame conditions and constancy of the domains imposed as premises. ($^{\ast}$)~In the \textbf{UF1} row, monotheism
holds in the qualified form: any two Godlike entities \emph{existing} at a
world are identical there.}
\label{tab:simplified}
\centering
\footnotesize
\setlength{\tabcolsep}{3.5pt}
\begin{tabular}{llccc}
\hline\noalign{\smallskip}
Postulate & Further axioms & $\Box\,\we x.\God x$ & Collapse & Monotheism \\
\noalign{\smallskip}\hline\noalign{\smallskip}
\textbf{U1}: $\Pos$ a modal ultrafilter & \textbf{A2}, \textbf{A3} & yes & \textbf{no} & yes \\
\textbf{F1}: $\HF_\God$ a modal filter & none & yes & \textbf{no} & \textbf{no} \\
\textbf{UF1}: $\HF_\God$ a modal ultrafilter & none & yes & \textbf{no} & yes$^{\ast}$ \\
\noalign{\smallskip}\hline
\end{tabular}
\end{table}

\subsection{Maximality, uniqueness, collapse --- and circularity}
\label{sec:dials}

Section~5.2 of \cite{OdifreddiGomes2026} treats these matters as a package; they
come apart. Under G\"odel/Scott, monotheism --- at most one Godlike being
--- follows from \textbf{A1} and
$\mathbf{df.}\God$.\footnote{\isaref{ScottVariant}{MT}.} In the theory of
Proposition~\ref{prop:hauptfilter} it fails, with a two-entity countermodel,
while maximality is untouched: $\God$ is still defined as having \emph{all}
positive properties. The examples thus establish specific
non-implications --- maximality does not imply uniqueness
(Proposition~\ref{prop:hauptfilter}), and neither maximality nor principal
generation implies the collapse (Section~\ref{sec:ufvariant}) --- rather
than full mutual independence of the three, which is not claimed here. Even so, these non-implications suffice for the point at issue. The claim of \cite[\S5.2]{OdifreddiGomes2026} that emending the argument yields a divine being falling short of Anselm's \enquote{that than which nothing greater can be conceived} rests on two steps, and both fail. Anderson's and Fitting's emendations do not give up principality at all (Table~\ref{tab:variants}); and where principality is given up, as under \textbf{F1}, maximality is untouched, since $\God$ is defined by having all positive properties.

One interpretive limit should be marked. \enquote{Maximality} here means
maximality relative to the chosen positivity predicate: $\God$ has every
property that \emph{this} $\Pos$ counts positive. That is not automatically
Anselm's unrestricted comparative notion of greatness, and a variant could
preserve the formal maximality while impoverishing what positivity ranges
over. The objection to \cite{OdifreddiGomes2026} is accordingly that the
inference from non-principality to lost maximality does not go through, not
that the formal notion settles the Anselmian question.

Finally, the circularity charge that \cite[\S5.1]{OdifreddiGomes2026} adopt from MacIntosh \cite{MacIntosh1991}: the possibility premise --- a Godlike being is \emph{possible} --- already amounts to the conclusion that one exists. A formal observation is
available here, though a limited one. Under \textbf{F1} one can derive that a
Godlike being actually, and even necessarily, exists, yet one \emph{cannot}
derive that a Godlike being \emph{possibly} exists: $\Diamond \we x.\,\God x$
has a countermodel in the weakest normal logic~$K$, and becomes derivable
once every world has some accessible world --- seriality, modal logic~D ---
and hence also under reflexivity~(T).\footnote{\isaref{SimpleVariantHF}{T3\_indep},
\isaref{SimpleVariantHF}{T3\_in\_D}, \isaref{SimpleVariantHF}{T3\_in\_T}.}
This turns on worlds without successors, where every $\Diamond$-formula
fails, and is to that extent an artefact of $K$; MacIntosh's diagnosis is
made for the settings in which the argument is normally run, where every
world has successors and the possibility premise is substantive. The observation should therefore be read as delimiting the formal scope of that analysis --- the equivalence of the possibility premise with the conclusion is not a logical triviality independent of the modal base, since in $K$ even necessary existence does not yield possible existence --- and not as an answer to it.
Whether the premise can be justified independently of the conclusion is a
separate question, which the present formal apparatus does not settle.

\section{Two Clarifications to Odifreddi and Gomes}
\label{sec:corrections}

\subsection{Theorem IV: a question of which premises, not a prover failure}
\label{sec:thmIV}

Section~2.3 of \cite{OdifreddiGomes2026} reports that, after the repair of
the inconsistency in G\"odel's original axioms, attempts to continue the
proof with formal verifiers \enquote{were not successful}: the automated
provers \enquote{fail to prove the crucial Theorem~IV} --- the possibility
of a Godlike entity, $\Diamond \we x.\,\God x$. The impression this leaves
--- that of an open failure of the tools --- is mistaken, but the
responsibility for it is not theirs alone: the wording echoes
\cite[\S4.4]{BenzmuellerScott2025}, co-authored by the present author, which
puts the matter in a way that invites exactly this reading. What follows is
therefore a clarification rather than a correction of an error.

Theorem~IV is not derivable when G\"odel's conjunction axiom
\textbf{Ax1} is restricted to its two-argument (\enquote{binary}) form --- and
this, too, is a theorem, not a prover failure. There is a countermodel in
which positivity is a \emph{non-principal ultrafilter} over countably many
individuals --- the natural numbers, all of which exist at the single world
of the model. That world accesses itself, so every frame condition from $K$
to $S5$ is satisfied and the result does not depend on the modal base; and
the domain is constant. In this model all axioms with binary \textbf{Ax1}
hold, yet the positive properties have empty intersection and no individual
is Godlike. The ultrafilter is obtained as a \emph{free} ultrafilter on the
natural numbers, taken from the Isabelle distribution, whose existence rests
on Zorn's lemma: the axiom of choice is the one non-constructive ingredient of the construction. This is an explicitly constructed model, not the output of a model finder.\footnote{The countermodel is constructed explicitly in
\thy{Th4Underivability}: \isaref{Th4Underivability}{Ax1}--\isaref{Th4Underivability}{Ax4}
verify the axioms in it, and \isaref{Th4Underivability}{Th4\_fails} that
Theorem~IV fails at every world; the ultrafilter comes from
\texttt{HOL-Nonstandard\_Analysis.Free\_Ultrafilter}.}

No finite countermodel exists --- where \emph{finite} must be read as
finitely many individuals \emph{and} finitely many worlds. What the argument
requires is that the space of property intensions
$\gamma = e \Rightarrow \sigma$ be finite, and over infinitely many worlds
that space is infinite even for a single individual; the
finite-ultrafilter observation of \cite[Thm.~3.2]{OdifreddiGomes2026} ---
every ultrafilter on a finite set is principal --- therefore does not
transpose to the modal setting as it stands. Whether the second hypothesis
can be dropped --- whether Theorem~IV already follows from the binary axioms
over finitely many individuals but infinitely many worlds --- is not settled
here; the infinite countermodel above has one world, so it does not decide
that case. Under both finiteness assumptions the observation does
transpose: the positive properties are then principally generated
and Theorem~IV follows, from the \emph{binary} \textbf{Ax1} together with
\textbf{Ax2a}, \textbf{Ax2b}, \textbf{Ax4} and
symmetry.\footnote{\isaref{GoedelAx1Gen}{Th4\_finite}, whose two hypotheses
are the finiteness of the individuals and of the worlds; it is stated and
proved before \textbf{Ax1Gen} is introduced in that theory, so it cannot appeal to it. This also delimits the reach of
finite-domain analyses generally: over a finite domain, principality ---
and with it much of the ultrafilter reading --- is automatic. On the
author's view the argument should ultimately be assessed over domains
comprising the infinitely many objects of mathematics, in keeping with
G\"odel's realism; see \cite{MuehlenbeckBenzmueller2026} and
\cite{BenzmuellerKirchnerPasetto2026} for this line.} This also explains the behaviour of the tools:
a model finder that enumerates finite structures, like Nitpick, can never exhibit such an infinite countermodel, and no finite one exists, while the provers cannot
derive the theorem, since the infinite countermodel exists; the question is
settled only by the explicit construction.

G\"odel's binary axioms therefore \emph{permit} precisely the
configuration that Friedman's construction exploits \cite{Friedman2012}: a
non-principal ultrafilter of positive properties whose total intersection
is empty. What rules it out is the footnote
G\"odel added to \textbf{Ax1} in his 1970 manuscript \cite{Goedel1970}
(\enquote{and for any number of summand[s]}), formalised as \textbf{Ax1Gen}:
\emph{every} conjunction of positive properties, including infinite ones, is
positive \cite{AndersonGettings1996,KanckosLethen2021}. From \textbf{Ax1Gen} the
lemma $\Pos(\God)$ --- Godlikeness is itself positive --- follows, and with \textbf{Ax2a} and \textbf{Ax4} a single automated step
proves Theorem~IV; the \emph{necessary} existence of a Godlike entity, $\lfloor\Box\we x.\,\God x\rfloor$ --- Theorem~V in the numbering of \cite{BenzmuellerScott2025}; \cite{OdifreddiGomes2026} list only Theorems~I--IV --- follows.\footnote{\isaref{GoedelAx1Gen}{L}; \isaref{GoedelAx1Gen}{Th4}; \isaref{GoedelAx1Gen}{Th5}.} Indeed,
\cite[\S2.3]{OdifreddiGomes2026} itself completes the proof along exactly
this route: \enquote{With this generalisation of Axiom~I and the definition
of a God-like entity, it becomes possible to prove [\dots] that being
God-like is a positive property. From this, it is possible to prove
Theorem~IV using Axioms~II and~IV} --- the same footnote generalisation,
the same lemma $\Pos(\God)$; and \enquote{Axioms~II and~IV} are
\textbf{Ax2a} and \textbf{Ax4}, the labels their own proof display uses as
well.

The clarification matters for the ultrafilter reading: \textbf{Ax1Gen} is
closure under arbitrary \emph{conjunctions} of positive properties --- under
the extensional reading, closure of the ultrafilter under arbitrary
\emph{intersections} --- without it the principality step
in \cite[\S3.1]{OdifreddiGomes2026}, resting on their Axiom~3.6
(\enquote{Being God is a positive property}), loses its derivation, and its
absence is what leaves Friedman's non-principal reading open.

\subsection{The extensionality of positivity}
\label{sec:ext}

To prove Theorem~IV, Odifreddi and Gomes \cite[\S2.3]{OdifreddiGomes2026}
invoke a premise they call \enquote{the extensionality of positivity
(derivable from Ax2a and Ax4)} --- properties with the same extension are
alike positive --- used once, to pass from the empty property
$\lambda x.\bot$ to the equally empty $\lambda x.\,\neg(x = x)$.
(\textbf{Ax2a}: of a property and its complement, exactly one is positive;
\textbf{Ax4}: a property necessarily entailed by a positive property is
positive.) Whether the premise holds depends on the reading of
\enquote{same extension}, and the two readings differ exactly for
non-rigid intensions: two concepts may pick out the same individuals at
the world of evaluation and yet come apart at another world --- think of
\emph{being the tallest person in the room} and \emph{being the oldest
person in the room} when both happen to apply to the same individual.
Extensionally this difference is invisible; intensionally it decides
whether positivity may separate the two concepts. Written out, the two
readings of \enquote{$X$ and $Y$ have the same extension} are
\[
  \text{(local)}\quad \wa z.\,(X z \leftrightarrow Y z)
  \qquad\text{versus}\qquad
  \text{(global)}\quad \Box\,\wa z.\,(X z \leftrightarrow Y z),
\]
and correspondingly two readings of the extensionality premise,
$\forall X Y.\,(\dots) \supset (\Pos X \leftrightarrow \Pos Y)$. Read
\emph{locally} the premise does not follow from \textbf{Ax2a} and
\textbf{Ax4}: Nitpick finds a two-world countermodel in
which two properties coincide at the world of evaluation, diverge at the
other world, and receive different positivity verdicts --- a countermodel that survives the $S5$ frame conditions with constant domains;\footnote{\isaref{ExtensionalityTests}{PosExt\_from\_Ax2a\_Ax4};
under the $S5$ frame conditions
\isaref{ExtensionalityTests}{PosExt\_from\_Ax2a\_Ax4\_S5}.} were the local
schema valid, the intensional reading of $\Pos$ would collapse into the
extensional one. Read \emph{globally} --- necessarily coextensive --- it
does follow, from \textbf{Ax4} alone: necessary coextension is mutual necessary
entailment, and \textbf{Ax4} applied in both directions transfers positivity;\footnote{\isaref{ExtensionalityTests}{PosExtNec\_from\_Ax2a\_Ax4}.} and only this
weaker schema is needed for their single application, since
$\lambda x.\bot$ and $\lambda x.\,x \neq x$ coincide at every world.\footnote{\isaref{ExtensionalityTests}{SelfDiff\_equiv\_Empty}.} So the parenthetical
is right on the global reading and wrong on the local one --- a distinction
the extensional idiom cannot even express. The appeal is in any case
dispensable: that the empty property is not positive follows from
\textbf{Ax2a} and \textbf{Ax4} directly, for the empty property necessarily
entails every property, in particular its own complement; by \textbf{Ax4}
its positivity would thus spread to its complement, contradicting
\textbf{Ax2a}.\footnote{\isaref{ExtensionalityTests}{NegEmpty}.}

Their proof of Theorem~IV is thus untouched; what is ambiguous is the general statement. The genuine demodalisation lies elsewhere --- in the general
ultrafilter interpretation, where G\"odel's intensional $\Pos$ gives way to
Friedman's extensional $\mathrm{POS}$ (explicit in their Section~4). That
move is legitimate for Friedman's purposes, but it is a change of system:
from a $\gamma$- to a $\delta$-ultrafilter, with immediate consequences for
the collapse.

The justification \cite[\S4]{OdifreddiGomes2026} offers for extensionality
repays the same analysis. The argument is Friedman's: if $X$ and $Y$ are
coextensive and $X$ is positive while $Y$ is not, then $X \sqcap {\sim}Y$
would be a \emph{positive empty property}. Friedman does not claim that the
axioms alone exclude such a property; he excludes it by the premise that
every positive property holds of some object, which he obtains from the
existence of God. That premise is modally innocent once the existence of a
Godlike individual at every world is in hand --- the Godlike individual
existing at $w$ has every property positive at $w$ --- and with it the local
step closes: the \emph{local} extensionality of positivity is a theorem of
G\"odel's full system, from binary \textbf{Ax1}, \textbf{Ax2a} and the
existence of God.\footnote{\isaref{GoedelAx1Gen}{PosExt\_local}, from
\isaref{GoedelAx1Gen}{ExistsG}.} So the argument is sound where Friedman
places it, after the existence proof. What local extensionality is not is a
consequence of \textbf{Ax2a} and \textbf{Ax4}, as the parenthesis of their
Section~2.3 --- read locally --- asserts; it lies downstream of the existence of God, that is, of
the very conclusion the proof of Theorem~IV is on its way to. That is why
only the global instance, available from \textbf{Ax4} alone, can be used at that point without circularity.

A related point: on the usual reading of necessary entailment, the empty
property entails \emph{every} property --- the root of the inconsistency
Leo-II detected in G\"odel's original axioms
\cite{BenzmuellerWP2016,BenzmuellerWP2016KI} via the \enquote{empty essence
lemma}; the falsum derivation is replayed in the ancillary files, in base
logic~$K$.\footnote{\isaref{GoedelInconsistency}{Inconsistency}.} A second repair
besides Scott's, which guards entailment against the empty property
\cite[\S4.5]{BenzmuellerScott2025}
(cf.\ \cite[\S3.1 and fn.~13]{FuenmayorBenzmueller2017}), blocks that lemma while
leaving G\"odel's own definition of essence untouched. Its consequences are
as follows:
the guarded theory is consistent and the falsum derivation disappears; at
the same time a degenerate one-element model appears in which the empty
property is positive and nothing is Godlike --- and it is once more
\textbf{Ax1Gen} that excludes precisely this
model.\footnote{Verified in \thy{ExtensionalityTests} and
\thy{AltEntailmentVariant}: \isaref{AltEntailmentVariant}{SecondModel};
\isaref{AltEntailmentVariant}{SecondModel\_excluded}.}

\section{The Friedman Connection, and Conclusion}
\label{sec:friedman}

None of the above touches Section~4 of \cite{OdifreddiGomes2026}, which is where their substantive new contribution lies. Nothing here verifies it: it has not been formalised, and the remarks below are invitations rather than objections.\footnote{One remark for the record. Their p.~107 displays $\mathrm{Con}(T_5) \leftrightarrow \mathrm{Con}(\mathrm{ZFC})$; the parenthesis that follows it notes that the right-to-left direction uses the existence of a measurable cardinal, so that what their results establish is the chain $\mathrm{Con}(\mathrm{ZFC} + \text{measurable}) \rightarrow \mathrm{Con}(T_5) \rightarrow \mathrm{Con}(\mathrm{ZFC})$, and the biconditional is to be read together with that parenthesis.} Three questions seem worth pursuing.

First, can Friedman's construction be easily formalised? His Divine Object axiom ranges only over the \emph{parameter-free definable} positive classes --- a countable
family whose intersection can be inhabited even when that of \emph{all}
positive classes is empty. Definability is a statement about the object
language itself, without direct counterpart in a shallow embedding; a
\emph{deep} embedding, which represents that object language explicitly,
could express it. Recent work provides deep and shallow embeddings of
first-order modal logic in HOL with machine-checked faithfulness
\cite{BenzmuellerKirchner2026}; extended to full higher-order modal logic,
it would offer a natural basis. Reconstructing Friedman's system~$T_5$ in
Isabelle/HOL is a particularly interesting open technical problem the paper
raises.

Second, how do the ways out compare? Friedman gives up \emph{principality}
(removing the perfect object makes the ultrafilter non-principal and thereby
mathematically powerful); Anderson and Fitting keep principality (Table~\ref{tab:variants}) and instead give up what forces the \emph{collapse}: the intensional reading of positivity in Fitting's case, the rigidity of the remaining ultrafilter $\Pos'$ in Anderson's. Are
these genuinely distinct manoeuvres, or are they closely related in some
way?

Third, where should the construction attach? It needs positivity read
extensionally, and Fitting's variant provides exactly that
extensionalisation, without collapse; starting there would also spare
\cite{OdifreddiGomes2026} the claim that collapse and circularity are
unavoidable (cf.\ also \cite{MuehlenbeckBenzmueller2026} on maximality and
collapse).

In conclusion, two theses must be kept apart. First, the historical thesis
--- that G\"odel himself did not regard the modal collapse as a defect,
defended by Kova\v{c} \cite{Kovac2012} and endorsed, among others, in
\cite[\S4.4]{BenzmuellerScott2025} and \cite{MuehlenbeckBenzmueller2026}
--- neither needs the ultrafilter
argument nor is damaged by its failure. Second, the structural thesis ---
that the collapse is inherent in \emph{any} ultrafilter characterisation of
positivity with God as principal generator --- is refuted by
Proposition~\ref{prop:ufvariant} and the explicit $S5$ model of
Section~\ref{sec:ufvariant}. That model carries the refutation on its own;
Proposition~\ref{prop:hauptfilter} and the \textbf{UF1} strengthening of
Table~\ref{tab:simplified} address the separate claims about principality
and monotheism, and Table~\ref{tab:variants} adds that the collapse-free
emendations of Anderson and Fitting keep principal generation by Godlikeness
as well. What drives the collapse is not the filter structure but the
rigidity axiom \textbf{A4} (the first conjunct of their Axiom~IIB) --- a
\emph{modal} condition,
as the opening observed --- and, by
Proposition~\ref{prop:equivalence}, under the ultrafilter postulate alone the two coincide: the collapse is the rigidity of positivity, with no help from G\"odel's essence and necessary-existence machinery. This comment thus takes up exactly the invitation of
\cite{OdifreddiGomes2026} --- to pursue the ultrafilter interpretation and
its surprising mathematical consequences --- if not quite in the
anticipated direction: the interpretation's central claim turns out to rest on a modal condition, rigidity, which its set-theoretic vocabulary does not express.

\section{Reproducibility, and the Status of the Machine Checks}
\label{sec:repro}

The companion session is \texttt{GoedelUltrafilterNote}: fifteen theory
files and a \texttt{ROOT}, built with \texttt{isabelle build -D . -v} under
Isabelle2025-2 in about a quarter of a minute. It depends on \texttt{HOL}
and, for the free ultrafilter of Section~\ref{sec:thmIV}, on
\texttt{HOL-Nonstandard\_Analysis}; Nitpick with its Kodkod backend and the
SMT solver used by a few proof steps ship with Isabelle. The build also
typesets the theories as a document. Exit code~$0$ certifies that every
proof closes and that every model-finder call meets the outcome recorded in
its \texttt{expect} annotation. The session accompanies this article as
ancillary files of \texttt{arXiv:2608.07578} (v4; \texttt{doi:10.48550/arXiv.2608.07578}), a permanent and versioned record
of the exact files on which the claims rest, and is offered to the journal
as electronic supplementary material.

\paragraph{What is checked, and how.} \enquote{Machine-checked} covers
three things here. (i)~\emph{Kernel-checked proofs}: $125$ of the $170$
statements are theorems accepted by the Isabelle/HOL kernel, among them
every positive result of this comment --- the derivabilities of
Propositions~\ref{prop:ufvariant} and~\ref{prop:hauptfilter}, the
equivalence of Proposition~\ref{prop:equivalence}, the principality entries
of Table~\ref{tab:variants}, and the extensionality facts of
Section~\ref{sec:ext}. (ii)~\emph{Explicitly constructed models, also
kernel-checked}: \thy{UFilterS5}, the $S5$ model of
Section~\ref{sec:ufvariant}, and \thy{Th4Underivability}, the infinite model
of Section~\ref{sec:thmIV}, interpret the types and constants of the
embedding and prove the axioms and the refutation in that interpretation,
so that the countermodel claim is itself a theorem. The first began as a
Nitpick countermodel and was then written out; the second could not, since
no finite countermodel exists. (iii)~\emph{Model-finder output}:
twenty-eight negative claims rest on Nitpick \cite{Nitpick} countermodels
and sixteen consistency claims on Nitpick satisfiability witnesses; one
check records that no model exists within a stated finite scope. Each call
fixes its scope --- the cardinalities of worlds and of individuals --- and
carries \texttt{user\_axioms}, without which \texttt{expect = genuine} would
be met vacuously by a model that ignores the theory's axioms (the two
exceptions are sanity checks in \thy{HOML}, which postulates nothing), and
Nitpick reports all of them as \emph{genuine}, not \emph{potentially
spurious}. They are finite, and a finite countermodel says nothing about
infinite models; any of them could be promoted to kind~(ii) by writing the
model out, as was done for the claim that carries this comment.

\paragraph{Witnesses and the object language.} Section~\ref{sec:a4} raised
the question whether a refutation uses a witness expressible in the modal
object language or one that names a world: the shallow embedding's
quantifiers range over domains that contain the nominals and satisfaction
operators of hybrid logic, so a proof may use resources the modal object
language lacks. This can be audited mechanically with the detector of
\cite[\S3]{BenzmuellerNominals2026}, which traverses the recorded proof
terms and reports, for each theorem, the terms substituted for variables of
lifted type, flagging those outside the grammar of the object language. For
the theories carrying the diagnosis of Section~\ref{sec:a4} the outcome is
this: in \thy{UFilterOnly} (Proposition~\ref{prop:equivalence}) and
\thy{UFilterVariant} (Proposition~\ref{prop:ufvariant} and the collapse
derivation from \textbf{A4}) no flagged witness occurs --- the proof of
Proposition~\ref{prop:equivalence} instantiates only the constant property
$\lambda x.\,\varphi$, schematic in $\varphi$; in \thy{UFilterS5} the
refutation at the atom uses only $q$ and $\lambda x.\,q$, whereas the
\enquote{fails at every world} readings use the world-indexed
$\lambda u.\,u = w$ and are flagged. Some world-naming is unavoidable when a
model is exhibited; the point is that which is which can be checked rather
than asserted.

\section*{Acknowledgements}
The author thanks David Fuenmayor, Bruno Woltzenlogel Paleo, Dana Scott
and Andrea Vestrucci for earlier collaborations in this line of work.

\bibliographystyle{amsplain}
\bibliography{note}

\end{document}